\documentclass[12pt,reqno]{amsart}
\usepackage{amsfonts}

\usepackage{amscd}
\usepackage{amsmath}
\usepackage{graphicx}

\theoremstyle{plain}

\theoremstyle{plain}
\newtheorem{Thm}{Theorem}

\newtheorem{ex}{Example}
\theoremstyle{definition}
\newtheorem{Def}{Definition}
\newtheorem{pf}{Proof}

\theoremstyle{remark}

\title{How to Distinguish a Local Semigroup from a Global Semigroup}  
\author{J. W. Neuberger}      
\date{6 September 2011} 
\email{jwn@unt.edu} 
\address{Department of Mathematics \\
University of North Texas \\
Denton, TX 76203-5017}     

\begin{document} 
\maketitle 
\centerline{\bf To Jerry Goldstein}
\medskip

\begin{abstract}
For a given autonomous time-dependent system that generates either a global, in time,
semigroup or else only a local, in time, semigroup, a test involving a linear eigenvalue
problem is given which determines which of `global' or `local' holds.  Numerical examples
are given.  A linear transformation $A$ is defined so that one has `global' or
`local' depending on whether $A$ does not or does have a positive
eigenvalue.  There is a possible application to Navier-Stokes problems.
\end{abstract}

\section{Introduction}\label{intro}
\begin{Def}

Suppose that $X$ is a complete separable metric space (i.e., a Polish space).
A function $T$ with domain $[0,\infty)$ and range in the collection of all transformations
from $X \rightarrow X$ is called a semigroup provided that 
\begin{itemize}

\item{$T(0)x = x,\; x\in X$ .}
\item{$T(t) T(s) = T(t+s), \; t,s \ge 0$ ($T(t)T(s)$ is the composition
of $T(t)$ and $T(s)$.)}
\item{$T$ jointly continuous, i.e., if 
\begin{equation*}
g(t,x) = T(t)x, \; t \ge 0, \; x \in X,
\end{equation*}
then $g$ is continuous.}
\end{itemize}
\end{Def}

Such a semigroup is called {\it linear} if $X$ is a linear space
and each $T(t), \; t \ge 0$, is a linear transformation.  Otherwise, $T$ is
called nonlinear.

\begin{Def}
Suppose $X$ is a subset of a Banach space $Y$, and $T$ is a semigroup on $X$.
The conventional generator of $T$ is
\begin{equation*}
B =  \{ (x,y) \in X \times Y  : \newline  y = \lim_{t \rightarrow 0+} \frac{1}{t}(T(t)x - x) .
\end{equation*}
\end{Def}

\section{Some History of One-Parameter Semigroups}\label{history}
  It has long been of interest to consider semigroups and their conventional
generators, specifically to recover a semigroup from its conventional generator.
Given a class of semigroups on a subset of a Banach space, one might try to
characterize the conventional generators of members of this class and then to
give a constructive means of recovering a semigroup in the class
from its conventional generator.
\medskip

A special case of the Hille-Yosida-Phillips theorem is the following:
\begin{Thm}
Suppose that

\begin{itemize}
\item{$X$ a Banach space.}
\item{$T: [0,\infty) \rightarrow L(X,X), \; t \ge 0.$}
\item{$T$ is a jointly continuous semigroup on $X$.}
\item{$|T(t)| \le 1, \; t \ge 0$ (nonexpansive).}
\item{$B$ is the conventional generator of $T$}.
\end{itemize}

Then
\begin{itemize}
\item{$B$ is a closed, densely
defined linear transformation on $X$,}
\item{$(I - \lambda B)^{-1} \in L(X,X)$ is nonexpansive, $\lambda \ge 0$, and}
\item{$\lim_{n \rightarrow \infty} (I - \frac{\lambda}{n} B)^{-n} x = T(\lambda)x, \hbox { for } 
 x \in X, \lambda \ge 0.$}
 \end{itemize}
 
 Moreover, if $B$ satisfies the first two items of the second group above, then
 
 \begin{itemize}
\item{there is a unique jointly continuous
nonexpansive linear semigroup $T$ whose conventional generator is $B$ (and hence
the last item above holds also).}
\end{itemize}

\end{Thm}
By the mid 1950s, a great deal was known about the relationship between linear
semigroups and their generators, (cf \cite{hp}) and, for a much more recent source, \cite{Goldstein}.  In the later 1950s, there was a fledging
effort to find analogies of the linear theory for nonlinear semigroups.  The books \cite{Brezis} and \cite{daprato} give an early 1970s
summary of how this development by analogy progressed.  Not so much
of a fundamental nature along the lines of analogy have transpired since.  Starting in the early 1970s, there were the first rumblings of an alternative theory, using a notion of generator that
can be traced to Gauss,
Riemann and Lie, with a resultant rather satisfactory theory finally given in \cite{DN}.  The main
result in \cite{DN} is given later in this note.
See \cite{probbook} for a much fuller historical account, with many references, on the development of linear and nonlinear semigroup theory.
\medskip

It was an early dream, dating from the late 1950s to eventually incorporate local semigroups
into a general framework.  Virtually nothing happened in this regard until recently.  The paper \cite{locsg} gives an account of some recent developments and so does \cite{probbook}.
The main focus in the present work is a recent update on \cite{locsg}, including some
 preliminary numerics which aim at the eventual gathering of numerical evidence to help
 decide whether a given semigroup is local or global.  The next section gives some
 thoughts on
 \cite{DN}, including relationships with Riemannian geometry.  Section \ref{locsalsg} gives a more detailed account of Lie generators for local semigroups and Section \ref{numerics} gives
 some numerical considerations.   

\section{Generators in the sense of Gauss, Riemann and Lie}\label{grlgen}

Suppose that $T$ is a jointly continuous semigroup on the Polish space $X$.  Denote by
$CB(X)$ the Banach space, under sup norm, of all bounded continuous real-valued 
functions on $X$.  Then
\begin{equation}\label{liegen}
A=\{(f,g) \in CB(X)^2 :
 g(x) = \lim_{t \rightarrow 0+} \frac{1}{t}(f(T(t)x) - f(x)), \; x \in X \}.
\end{equation} 
is called the Lie generator of $T$ (more fully, the Gauss-Riemann-Lie generator).
Why attach these names to this generator?  For Sophus Lie it come directly from
his work \cite{slie} although not for the purposes of the present paper.  He was
trying to establish a theory of ordinary differential equations by means of integrating
factors.  His generators were an essential part of his constructions.   But why the names of Gauss
and Riemann?  In their constructions of new geometries, they devised how
 to use a differentiable structure on what were essentially topological spaces.  It seems
 clear that their work inspired Lie's work.  How does
 one differentiate a real-valued function $g$ on a Riemannian manifold $M$?  One may take a function
 $f: [-1,1] \rightarrow M$ and then attempt to differentiate the composition $f(g)$.  This composition
 is a real-valued function on $[-1,1]$.  It at least make sense to ask whether or not it has
 a derivative.  A main point of divergence that the present work makes with 
 Riemannian geometry, is that in Riemannian geometries,  one has the notion of a differential function $g$ on $M$, and
 a differentiable function $f$ from $[-1,1] \rightarrow M$.  These are given in terms of charts and atlases.  In this context, such a composition $f(g)$
 is automatically differentiable.  In the case of \eqref{liegen}, a composition
 \begin{equation*}
 f(T(\cdot)x)
 \end{equation*}
 there may or may not be a bounded continuous function $g$ so that for all $x \in X$,
 \begin{equation*}
 (f(T(\cdot)x)^{\prime}(0) = g(x)
 \end{equation*}  If this does happen, then $f$ is admitted to the domain of $A$, otherwise not.  A key
 point is that in defining $A$ in \eqref{liegen},  one uses only differentiation that is earned,
 not hypothesized.  This idea is very much in line with the second half of Hilbert's Fifth Problem.
 \medskip
 
The theorem to follow gives a characterization of Lie generators of jointly continuous semigroups.
It is needed for the argument in the following section, which deals with local semigroups.  It
was found by J. R. Dorroh and the present writer.

 \begin{Def}
A sequence $\{f_n\}^{\infty}_{n=1}$ in $CB(X)$  $\beta-$converges to $f \in CB(X)$ if  it
is uniformly bounded and converges uniformly to $f$ on compact subsets of $X$. 
 \end{Def}
 See \cite{DN} for details on this kind of convergence and its resulting topology.
 
 \begin{Thm}\label{dnthm}
 Suppose that $X$ is a Polish space, $T$ is a jointly continuous semigroup on $X$ and $A$
 is its generator in the sense of \eqref{liegen}.  Then
 
 \begin{itemize}
\item{$A$ is a derivation.}
\item{$D(A)$ is $\beta-$dense in $CB(X)$.}
\item{If $\lambda \ge 0, (I - \lambda A)^{-1} \in L(CB(X))$ and is nonexpansive.}
\item {If $\gamma > 0$ then $ \{ (I - \frac{\lambda}{n} A)^{-n} : 0  \le \lambda \le \gamma, n=1,2,\dots \} $  is uniformly $\beta-$equicontinuous.}    
\end{itemize}
Moreover, if $A$ satisfies the four items above, then there is a unique jointly continuos semigroup
on $X$ which has $A$ as its Lie generator.
\end{Thm}
In a sense, this satisfies the quest for a complete generator-semigroup theory for the
class of jointly continuos semigroups on $X$:  Given such a semigroup $T$, there is a unique
generator, specified in terms of the above four items.  Given $A$ satisfying these four properties,
there is a unique $T$ which has Lie generator $A$.  See \cite{DN} for an argument.

\section{Local Semigroups}\label{locsalsg}
\begin{Def} 
The statement that $T$ is a local semigroup on the Polish space  $X$ means that
 $T$ is a function from a connected subset of
$[0,\infty)$ into the collection of all functions on $X$ to $X$ so that
the following hold:  
\begin{itemize}
\item{There is $m$,  a function from $X$ to $(0,\infty]$ such that $\frac{1}{m}$ is continuous
and is not identically equal to infinity,}
\item{$x \in  D(T(t)) \iff  t \in [0,m(x)) $, } 
\item{$\hbox{ If } t,s \ge 0, x \in X, \hbox{ then } T(t)T(s)x = T(t+s)x \iff t+s < m(x)$,}
\item{$T$ is jointly continuous and maximal ($\lim_{t \rightarrow s-} T(t)x$ exists
$\implies s < m(x)$)}
\end{itemize}
\end{Def}
A Lie generator $A$ for such a $T$ is defined as in \eqref{liegen}.
\medskip

  In \cite{locsg} there is a theorem which give a partial characterization
of the Lie generator of a local semigroup, that is, it gives some properties of such
a generator.    Developments there
can be used to give a proof of Theorem \ref{loc} to follow, but the argument given
here is much shorter.


\begin{Thm}\label{loc}
Suppose $T$ is either a local or a global jointly continuous semigroup  on $X$ and
$A$ is the Lie generator of $T$.
Then $A$ has a positive eigenvalue if and only if $T$ is local.

\end{Thm}

\begin{pf}
Suppose that $T$ is a local semigroup and define $f \in CB(X)$ by
\begin{equation*}
f(x) = \exp(-m(x)), \; x \in X.
\end{equation*}
Then for some $x \in X$ and $t \ge 0$,
\begin{equation*}
f(T(t)x) = \exp(-m(T(t)x)) = \exp(-(m(x) -t)) = f(x) \exp(t),
\end{equation*}
and so
\begin{align*}
& \frac{1}{t}(f(T(t)x) - f(x)) =  \\ & \frac{1}{t}(f(x) (\exp(t) - 1) \rightarrow f(x) = \exp(-m(x)) \hbox{ as } t \rightarrow 0+.
\end{align*}
Thus, 
\begin{equation*}
f \in D(A) \hbox{ and } (Af)(x) = f(x).
\end{equation*}
Since $f(x) \ne 0$ for at least one $x \in X$, it follows that $f$ is an
eigenvector of $A$ with eigenvalue one.
\medskip 

Now suppose that $A$ has a positive eigenvalue $\lambda$ with eigenfunction
$g \in CB(X)$, but that
$T$ is a semigroup (i.e., is not a local semigroup).  From Theorem \ref{dnthm},
it follows that
\begin{equation*}
(I - \lambda A)^{-1} 
\end{equation*}
exists.  But this is impossible since 
\begin{equation*}
(I  - \lambda A)g = 0
\end{equation*} 
and hence $g = 0$, a contradiction since the zero function is never an eigenfunction.
\end{pf}

\begin{Thm}\label{choice}
Suppose that $T$ is a local semigroup on $X$ with Lie generator $A$ and $\lambda >0$.
Then $\lambda$ is an eigenvalue of $A$.
\end{Thm}

\begin{pf}
As in Theorem \ref{loc}, if one starts with
\begin{equation*}
g(x) = \exp(- \lambda m(x)), \; x \in X, m(x) \ne 0,
\end{equation*}
one may conclude that
\begin{equation*}
Ag  = \lambda g.
\end{equation*}
\end{pf}

It seems mysterious that the Lie generator of a local semigroup has 
all of $(0,\infty)$ in its spectrum.
\medskip

\section{Three Examples}\label{examples}

\begin{ex}\label{example1}

\medskip

Suppose 
\begin{equation*}
X = [0,\infty) \hbox{ and }B(x) = x^2, x \in X.
\end{equation*}\label{de1}
For $x \ge 0$, the solution $u$ to
\begin{equation}
u(0) = x, \; u^{\prime} (t) = B(u(t)), \; t \in [0,\frac{1}{x}) 
\end{equation}
is given by
\begin{equation*}
 u(t) = \frac{x}{1 - t x}, x \ge 0, \; t \in [0,m(x))
 \end{equation*}
 where
 \begin{equation*}
 m(x) = \frac{1}{x}, \; x \in [0,\infty).
\end{equation*}
 One can see this by solving, for $x \ge 0$, the
 problem
 \begin{equation}\label{ex1}
 u(0) = x, \; u^{\prime} (t) = u(t)^2, \; t \in [0,m(x))
 \end{equation}
 and also noticing that $[0,m(x))$ is the maximal connected
 subset of $[0,\infty)$ containing zero over which there is a solution.
 Note that the local semigroup $T$ so that
 \begin{equation*}
 T(t)x = \frac{x}{1 - t x}, \; x \in [0,\infty), \; t \in [0,m(x)),
 \end{equation*}
 is associated with \eqref{ex1} in the sense that if $u$ satisfies
 \eqref{ex1}, then
 \begin{equation*}
 T(t) x =u(t), \; x \in [0,\infty), \; t \in [0,m(x)),
\end{equation*}
Thus defining $f \in CB(X)$ by
 \begin{equation}\label{mev}
f(x) = \exp(-\frac{1}{x}) = \exp(-m(x)), \; x \ge 0, \; f(0) = 0,
\end{equation}
$f$ yields an eigenvector of $A$ with positive eigenvalue.
Hence it is consistent with Theorem \eqref{choice} that $T$ is a local semigroup,
which may readily be seen anyway.
\medskip

There is a second way to connect $A$ with an eigenfunction $h$ of $A$ with eigenvalue one:
For such an $h$,
\begin{equation*}
h(x) = (Ah)(x) = \lim_{t \rightarrow} \frac{1}{t}(h(T(t)x) - h(x)) = h^{\prime}(x)B(x), \; x > 0,
\end{equation*}
so that 
\begin{equation}\label{heq}
 (A h)(x) = h^{\prime}(x) x^2, \; x > 0.
\end{equation}
If $h(x)$ is positive at some $x > 0$, then by
directly solving \eqref{heq} for $h$,
\begin{equation*}
h(x) = c \exp(-\frac{1}{x}), \; x > 0,
\end{equation*}
for some positive number $c$.  Clearly this agrees with
the result \eqref{mev} since by continuity,
\begin{equation*}
h(0)= \lim_{t \rightarrow} c\exp(-\frac{1}{x}) = 0.
\end{equation*}
Note that \eqref{heq} is a singular equation with the property that a constant times a 
solution is also a solution.
\medskip
\end{ex}

\medskip
\newpage
\begin{ex}\label{ex3}
Suppose 
\begin{equation*}
X = [0,\infty) \hbox{ and }B(x) = x(x-1), x \in X.
\end{equation*}
For $x \ge 0$, the solution $u$ to
\begin{equation}\label{de3}
u(0) = x, \; u^{\prime} (t) = B(u(t)), \; t  \ge 0) 
\end{equation}
is given by
\begin{equation*}
u(t) = \frac{x}{x + \exp(t) (1-x)}, \; t  \in [0,m(x))
\end{equation*}
where
\begin{equation*}
m(x) = 
\begin{cases}
\infty & \text{ if $ x \in [0,1],$} \\
\ln(\frac{x}{x - 1}) &\text{ of $ x >1.$}
\end{cases} 
\end{equation*}
The corresponding local semigroup $T$ is given by
\begin{equation*}
T(t)x =\frac{x}{x + \exp(t) (1-x)}, \; t  \in [0,m(x))
\end{equation*}
For $A$ the Lie generator of $T$, an eigenfunction for
$A$ with one as its eigenvalue is given,
following Theorem \ref{choice}, by
\begin{equation*}
f(x) = \exp(-m(x)) = 
\begin{cases}
0 & \text{ if  $x \in [0,1],$} \\
\frac{x-1}{x} &\text{ if $ x > 1.$}
\end{cases}
\end{equation*}
Note that any non-zero multiple of $f$ is also an eigenfunction
of $A$ for eigenvalue one.
\medskip

One can compute an eigenvector for this case following the same
path as the second approach of the previous example.
\end{ex}

\begin{ex}\label{ex2}
Suppose 
\begin{equation*}
X = [0,\infty) \hbox{ and }B(x) = -x^2, x \in X.
\end{equation*}
For $x \ge 0$, the solution $u$ to
\begin{equation}\label{de1}
u(0) = x, \; u^{\prime} (t) = B(u(t)), \; t  \ge 0) 
\end{equation}
is given by
\begin{equation*}
 u(t) = \frac{x}{1 + t x}, x \ge 0, \; t \ge 0.
 \end{equation*}
 Hence the corresponding semigroup (it is a semigroup,
 not a local semigroiup) is given by $T$:
 \begin{equation*}
 T(t)x = \frac{x}{1 + t x}, \; t,x \ge 0.
 \end{equation*}
 Denote by $A$ the Lie generator of $T$.
 If $A$ were to have a eigenfunction $h$ with eigenvalue one, say,
 then it would have to be that
 \begin{equation*}
 (Ah)(x) = f^{\prime}(x) (-x^2) = f(x), x > 0.
 \end{equation*}
 But solving
 \begin{equation*}
  f^{\prime}(x) (-x^2) = f(x), x > 0
  \end{equation*}
  for a non-zero $f$ which is positive somewhere on $(0,\infty)$ yields that for some $c > 0$,
  \begin{equation*}
  f(x) = c \exp(\frac{1}{x}), \; x > 0,
  \end{equation*}
 clearly unbounded with no hope of extending by continuity to $[0,\infty)$.
  This is also in agreement with Theorem \ref{loc}.
  \medskip
  
   \end{ex}
  The next section concerns some numerics and some comments on applications.
 
\section{Some Numerics}\label{numerics}
In this section there are numerical indications of eigenvectors for
finite dimensional approximations to each of the examples in the preceding
section.  In each case we truncate $X = [0,\infty)$ into an interval $[0,z]$ for
some $z > 0$ and divide this interval into $n$ pieces of equal length.  In each case, denote
by $A_n$ an appropriate finite dimensional approximation to the relevant Lie generator on the $n+1$ 
dimensional space $X_n$ of functions on $G_n = z*\frac{0}{n},z*\frac{1}{n} ,\dots,z*\frac{n}{n}$.
Define $\phi_n: X_n \rightarrow R$ by
\begin{equation*}
\phi_n (g) = \frac{1}{2} \| g - A_n g\|^2_{R^{n+1}}, \; g \in X_n.
\end{equation*}
It is straightforward enough to calculate an
ordinary gradient $\alpha_n(g)$ for $\phi_n$ evaluated at
an element $g \in X_n$.  Our method is one of steepest
descent to find a zero of $\phi_n$, but since $A$ is essentially a differential
operator on a function space, experience indicates that steepest descent
with this ordinary gradient is best avoided.  We instead use an appropriate
weighted Sobolev gradient $\nabla \phi_n$ which is effectively a preconditioned
version of $\alpha_n$ (see \cite{sv2}) for an account of Sobolev gradients).  To
get $\phi_n$ from $\alpha_n$ for a given one of Examples 1,2,3, first define $v_n$
from the appropriate function $B$:
\begin{equation*}
v_n (j) = B(\frac{zj}{n}), \; j=0,1,\dots, n.
\end{equation*}
Let $D$ be the linear transformation which takes first differences of members of $X_n$.
Define 
\begin{equation*}
Q: X_n \rightarrow X_n:  \; Q = I + (v*D)^t (v*D).
\end{equation*}
Finally define, for $g \in X_n$,
\begin{equation*}
(\nabla \phi_n)(g) = Q^{-1} ((\nabla \alpha_n)(g)).
\end{equation*}
The gradient $\nabla \phi_n$ is then used in steepest descent with optimal
step size to find a zero of $\phi_n$, that is an eigenvector of $A_n$ if the
limiting value of this iteration is not zero.
\medskip

Results for such an iteration in the case of each of the three examples are given below.
In the first two cases, the results track closely the explicit expressions found for eigenfunctions,
remembering that any non-zero multiple of an eigenfunction of $A$ is still an eigenfunction of $A$.
\medskip

In Figure 1, careful examination of the plot shows an indication of the zero derivative
at zero of the eigenfunction,  In Figure 2, the portion of the eigenfunction over $[0,1]$
that is zero is clearly indicated.  This indication was first noticed numerically, leading us
to realize that in Example 2, trajectories starting at $x \in [0,1]$ are global, that is
$m(x) = \infty$ for such $x$.  The indication that $T$ from Example 2 is a local semigroup
is reflected in the boundedness of the graph in Figure 2.   Extending the indicated
interval $[0,10]$ to a series of progressively longer intervals confirms this phenomena.
\medskip

Figure 3, for Example 3, needs more explanation.   A key to properly interpreting
Figure 3 is to notice that factor $10^{-6}$ at the top of the vertical axis.  The
graph thus indicates a numerical analyst's version of zero, so that Figure 3 does
not indicate an eigenfunction with positive eigenvalue but rather just the zero function,
which is no eigenfunction at all.
\medskip

These three results are supported by a large number of computer runs with different
interval lengths, and different number of grid points $n+1$.  Compelling numerical results are of
course not proofs of themselves, but rather strong pointers in the direction of proofs.  A further
discussions of relevant matters is in the next section.

\begin{figure}
\vspace{-3in}
\includegraphics[width=5.5in]{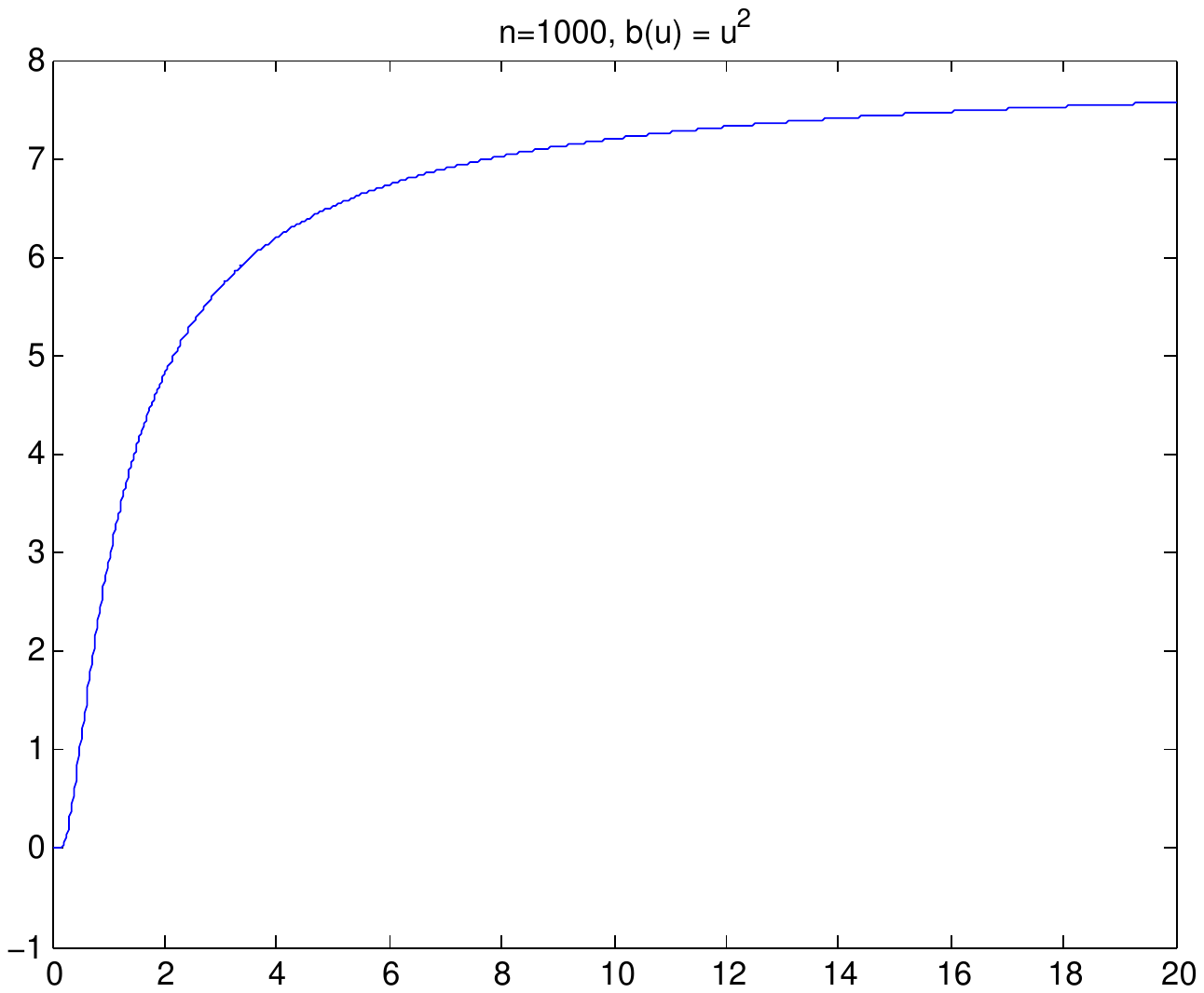}
\vspace{-2.in}
\caption{$u' = u^2$}
\end{figure}

\vspace{.5in}
\begin{figure}
\vspace{-1in}
\hspace{2.in}
\includegraphics[width=4.8in]{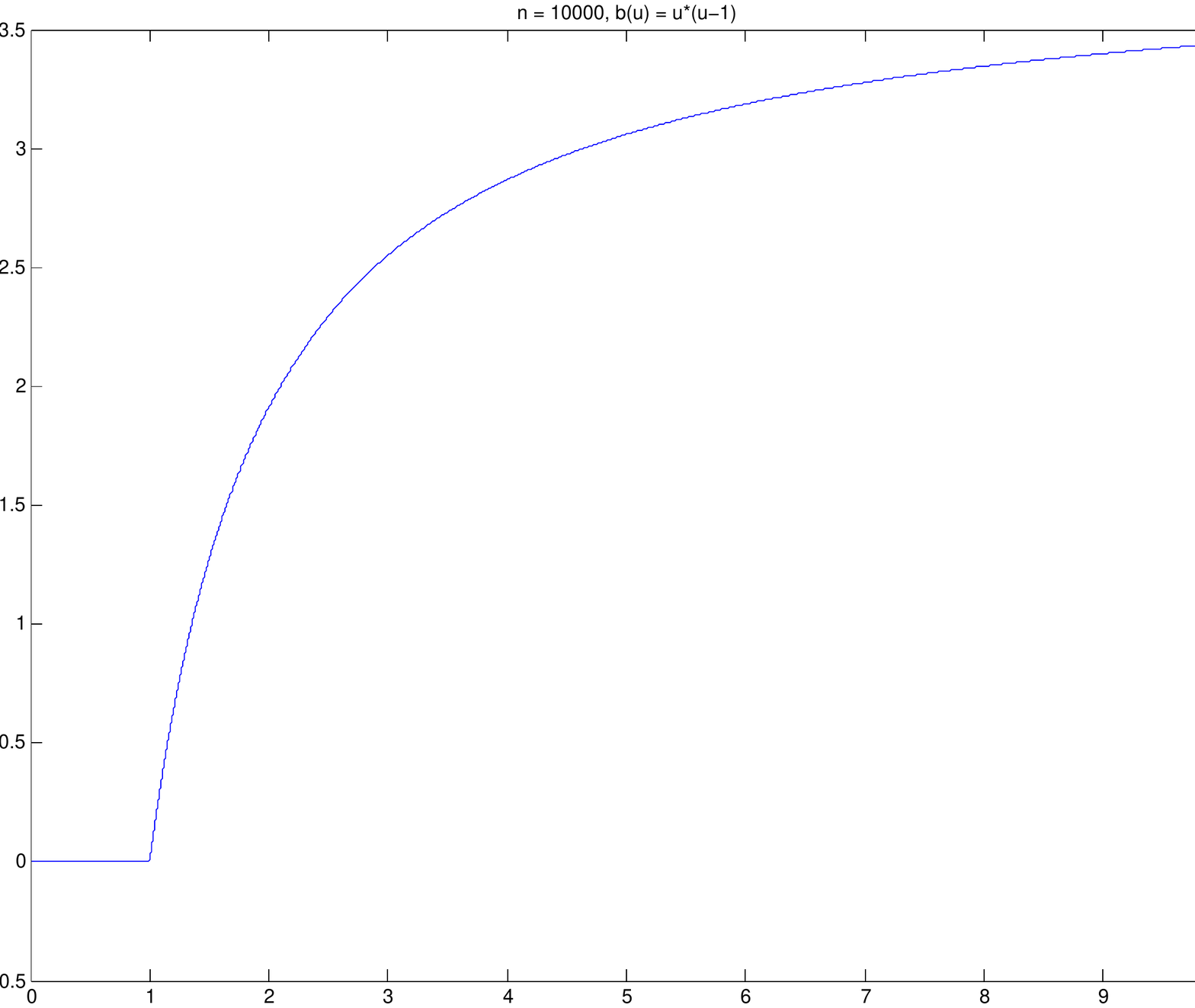}
\vspace{-1in}
\caption{$u' = u*(u-1)$}
\end{figure} 

\begin{figure}
\vspace{-1in}
\hspace{2.in}
\includegraphics[width=4.8in]{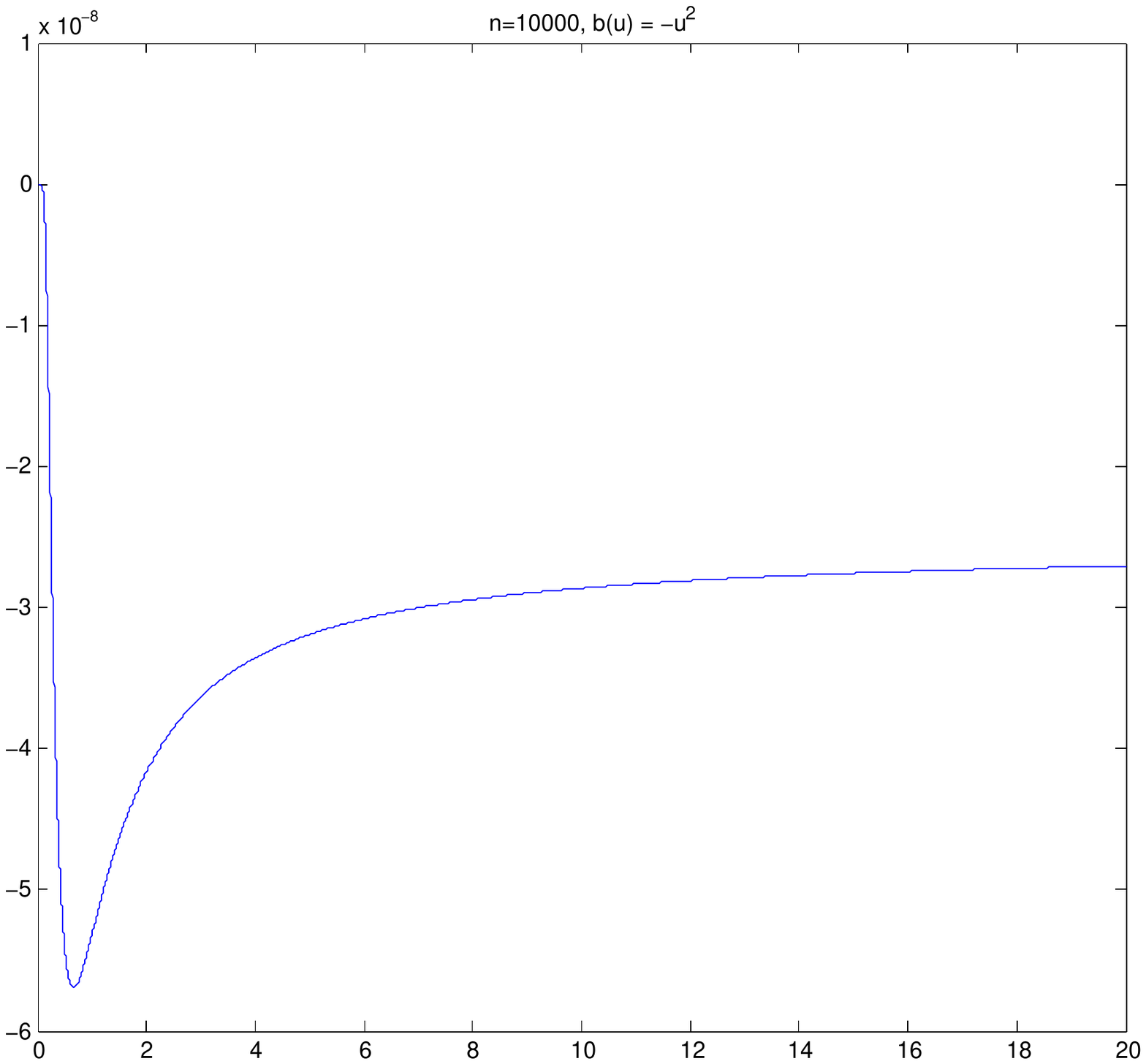}
\vspace{-1in}
\caption{$u' = -u^2$}
\end{figure} 




\newpage
\section{Epilogue}\label{epiloge}
What, if any, are some possible scientific and mathematical relevances of the above
development?  In Sophus Lie's book \cite{slie}, Chapter 4, there is the rather bold
statement that all information about a semigroup or local semigroup (my language -
Sophus Lie said simply `group') is contained  in its generator.  By this he clearly means
generators $A$ which share crucial properties with the `Lie generators' of the 
present paper.  The main point of the present
paper is to try enable some practical applications - here
the distinguishing of semigroups from local semigroups.  The `old' semigroup theory
starting from possibly \cite{N1} is summerized in \cite{Brezis} and \cite{daprato}.  It 
is rather limited in scope - being restricted to mainly nonexpansive semigroups (and definitely not local
semigroups) on a Hilbert space.  The 
semigroup theory as indicated in this paper does not have these restrictions.  
The example of Webb, \cite{webb}, demonstrated that a conventional generator,
even in very simple appearing semigroups, may have a domain so sparse that further 
analysis using this generator is discouraged.  Lie generators do not suffer from this.
\medskip

The ultimate evaluator of a semigroup theory is its utility in scientific problems as well
as its intrinsic mathematical value.  The present work seeks to indicate some of the former
concerns, the latter being left to individual tastes.
\medskip

Here we elaborate a bit on how the `old' and `new' theory may relate.  Semigroup
theory has always been an abstraction of time-dependent autonomous process, i.e.
processes evolving in time for which changes in a present state depend only on that
state and not what time it might be.  In an autonomous process the {\em law} of evolution is
not changing with time (no Federal Reserve changing interest rates on a bond which has continuously
compounding interest, in the middle of
the tenure of a bond, for example).  An autonomous process is often expressed as a
differential equation on some subset of a Banach space $X$:
\begin{equation}\label{geneq}
u(0) = x \in X, \; u^{\prime}(t) = B(u(t)), \; t \hbox{ in some interval containing zero}.
\end{equation}
Even systems as complicated as Navier-Stokes (with no external forcing term) may 
be expressed in this way.    What is a connection between $B$ in \eqref{geneq} and
a Lie generator of the semigroup or local semigroup $T$ connected with \eqref{geneq}?  It is this:
\begin{equation}\label{connect}
(Af)(x) = \lim_{t \rightarrow 0+} \frac{1}{t} (f(T(t)x) - f(x)), \; x \in X.
\end{equation}
Now in the presence of sufficient regularity of $f \in CB(X)$ and of $T$ (differentiability properties),
\eqref{connect} reduces to 
\begin{equation*}
(Af)(x) = f^{\prime}(x) Bx, \; x \hbox{ in the domain of $B$},
\end{equation*}
a relation that would have been instantly recognized by Sophus Lie.
Thus for a problem given by \eqref{geneq}, one has a rather concrete expression
for the Lie generator of its corresponding semigroup or local semigroup.  A discretization,
as in the examples, is not far behind, and thus numerics seeking to determine the
existence or not of a positive eigenvalue of $A$ are within reach.  Beyond problems
with a single or just a few ordinary differential equations given by \eqref{geneq}, 
appropriate discretizations involving Lie generators lead to problems of great size, quickly becoming a fatal strain on
even the largest of present day computers.  However, it is a productive
to code what can be run on present day machines, then lying in wait for bigger,
and particularly more efficient machines (in regard to communication between
processors)  to come along.
\medskip

Someone suspecting that significant bits of optimism expressed above are warranted,
might embark on a series of problems of increasing complexity.  First systems of
several ordinary differential equations (this is already started) and then time-dependent partial differential equations in one space dimension, then two, then three dimensions..
\medskip

We close with two general comments.  For dynamical systems and semidynamical systems
$T$ on a space $X$, many issues involve fixed points of $T$, i.e., points $x \in X$ so that
\begin{equation*}
T(t) x = x, \; t \ge 0 \hbox{ or } t \in [0,m(x)).
\end{equation*}
Examination of Examples 1-3 reveal that such fixed points of the corresponding $T$ play a
significant role when analysis is done with a Lie generator.  I suspect that this will be
a recurring item as use of Lie generators expands.
\medskip

A final comment has to do with the function $m$ and the space $CB(X)$ in the definition of
a local semigroup.  That $\frac{1}{m}$ is required to be continuous may be too strong.  A
weaker, unknown to me, condition on $m$, together with a relaxation of the requirement
that members of $CB(X)$ be continuous (a suggestion of someone attending my seminar at
the University of Texas a few years ago).  The suggestion that
$\frac{1}{m}$ be merely upper-semicontinuous  was made.  Possible application to
Navier-Stokes problems, for example, may depend on such extensions, the issue being
whether for such a system the continuity property on $\frac{1}{m}$ holds.  Accomodation
such an $m$ in an extension of the present setting might require an extension of $CB(X)$ beyond
continuous functions.  I hasten to add that I do not know that such an extension is needed
in order to treat Navier-Stokes using propositions of this paper.

\end{document}